\documentclass[12pt]{amsart}
\usepackage{amsmath, amssymb, latexsym}

\def\Alg{{\mbox{\rm Alg}}}
\def\End{{\mbox{\rm End}}}
\def\lcm{{\mbox{\rm lcm}}}
\def\Tr{{\mbox{\rm Tr}}}

\def\a{{\alpha}}
\def\b{{\beta}}

\def\w{{\omega}}

\def\BZ{{\mathbb Z}}

\def\tr{{\mbox{\rm Tr}}}

\numberwithin{equation}{section}

\newtheorem{theorem}{Theorem}[section]
\newtheorem{corollary}[theorem]{Corollary}

\newtheorem{proposition}[theorem]{Proposition}

\newtheorem{lemma}[theorem]{Lemma}

\begin{document}
\author{Siu-Hung Ng}
\title{Hopf Algebras of Dimension $2p$}
\address{\footnote{Current address}Department of Mathematics\\
        Iowa State University\\
        Ames, IA 50011 \newline
        Mathematics Department\\
        Towson University\\
        8000 York Road, Baltimore, MD 21252}
   \email{rng@math.iastate.edu,  rng@towson.edu}
\begin{abstract}
Let $H$ be a finite-dimensional Hopf algebra  over an
algebraically closed field of characteristic 0.
If $H$ is not semisimple and $\dim(H)=2n$
for some odd integer $n$, then $H$ or $H^*$ is not unimodular. Using this result,
we prove that if $\dim(H)=2p$ for some odd prime $p$, then
$H$ is semisimple. This completes the classification of Hopf
algebras of dimension $2p$.
\end{abstract}
\date{}
 \maketitle

\noindent
\section{Introduction}
In recent years, there has been some progress on the
classification problems of finite-dimensional Hopf algebras over
an algebraically closed field $k$ of characteristic 0 (cf.
\cite{Mont98}, \cite{And_bk2002}). It is shown in \cite{Zhu94}
that Hopf algebras of dimension $p$, where $p$ is a prime, are
isomorphic to the group algebra $k[\BZ_p]$. In \cite{Ng02} and
\cite{Ng02a}, the author completed the classification of Hopf
algebras of dimension $p^2$, which started in \cite{AS98} and
\cite{Mas96}. They are group algebras and Taft algebras  of
dimension $p^2$ (cf. \cite{Taft71}). However, the classification of Hopf algebras of
dimension $pq$, where $p, q$ are distinct prime numbers, remains open in
general.\\

It is shown in \cite{EG99}, \cite{GW00} that semisimple Hopf
algebras over $k$ of dimension $pq$ are trivial (i.e. isomorphic
to either group algebras or  the dual of group algebras). Most
recently, Etingof and Gelaki proved that if $p$, $q$ are odd prime
such that $p < q \le  2p+1$, then any Hopf algebra over $k$ of
dimension $pq$ is semisimple \cite{EG03}. Meanwhile, the author
proved the same result, using different method, for the case that $p$,
$q$ are twin primes \cite{Ng03}. In addition to that Williams
settled the case of dimensions 6 and 10 in \cite{WI}, and Beattie
and Dascalescu did dimensions 14, 65
in \cite{BD02}. Hopf algebras  of dimensions 6, 10, 14 and 65 are semisimple and
so they are trivial.\\

In this paper, we prove that any Hopf algebra of dimension $2p$, where $p$ is an odd prime, over
an algebraically closed field $k$ of characteristic 0, is semisimple. By \cite{Masuoka95}, semisimple Hopf
algebras of dimension $2p$ are isomorphic to
$$
k[\BZ_{2p}],\quad k[D_{2p}]\quad\mbox{or}\quad k[D_{2p}]^*
$$
where $D_{2p}$  is the dihedral group of order $2p$. Hence, our
main result Theorem \ref{main} completes the classification of Hopf algebras of
dimension $2p$.

\section{Notation and Preliminaries}\label{s1}
Throughout this paper, $p$ is an odd prime, $k$ denotes an algebraically closed
 field of characteristic 0, and $H$ denotes a finite-dimensional Hopf algebra over
$k$ with antipode $S$. Its comultiplication and counit are
respectively denoted by $\Delta$ and  $\epsilon$. We will use
Sweedler's notation \cite{Sw69}:
$$
\Delta(x) = \sum x_{(1)} \otimes x_{(2)}\,.
$$

A non-zero element $a \in H$ is called group-like if $\Delta(a)=a
\otimes a$.  The set of all group-like elements $G(H)$ of $H$ is a linearly
independent set, and it forms a group under the multiplication of
$H$. For the details of elementary aspects for
finite-dimensional Hopf algebras,  readers are referred to
the references \cite{Sw69} and \cite{Mont93bk}.\\

Let $\lambda \in H^*$ be a non-zero right integral of $H^*$ and  $\Lambda\in H$  a non-zero
 left integral of $H$. There exists $\a \in \Alg(H,k)=G(H^*)$, independent of
the choice of $\Lambda$,  such that $\Lambda a = \a(a) \Lambda$ for  $a \in H$.
Likewise, there is a group-like element $g \in H$, independent of
the choice of $\lambda$,  such that $\b\lambda  = \b(g)\lambda$ for  $\b\in H^*$.
We call $g$ the distinguished group-like element of $H$ and $\a$
the distinguished group-like element of $H^*$. Then we have Radford's
formula~\cite{Radf76} for $S^4$ :
\begin{equation}\label{eS}
S^4(a) = g(\a \rightharpoonup a \leftharpoonup \a^{-1}) g^{-1}
\quad\mbox{for }a \in H\,,
\end{equation}
where $\rightharpoonup$ and $\leftharpoonup$ denote the natural
actions of the Hopf algebra $H^*$ on $H$ described by
$$
\b\rightharpoonup a = \sum a_{(1)}\b(a_{(2)}) \quad \mbox{and}
\quad a \leftharpoonup \b = \sum \b(a_{(1)})a_{(2)}
$$
for $\b \in H^*$ and $a \in H$.  In particular, we have the following proposition.
\begin{proposition}[\cite{Radf76}]\label{p0}
  Let $H$ be a finite-dimensional Hopf algebra with antipode $S$ over the field $k$. Suppose that
  $g$ and $\a$ are distinguished group-like elements of $H$ and $H^*$ respectively. Then the order of $S^4$
  divides the least common multiple of the order of $g$ and the order of $\a$. \qed
\end{proposition}
For any $a  \in H$, the linear operator $r(a) \in \End_k(H)$ is defined by $r(a)(b)=ba$ for $b \in H$.
The semisimplicity of a finite-dimensional Hopf algebra can be characterized by the antipode.
\begin{theorem}[\cite{LaRa87}, \cite{LaRa88},  \cite{Radf94}] \label{t2}
  Let $H$ be a finite-dimensional Hopf algebra with antipode $S$ over the field $k$. Then the following
  statements are equivalent:
  \begin{enumerate}
    \item $H$ is not semisimple;
    \item $H^*$ is not semisimple;
    \item $S^2 \ne id_H$;
    \item $\tr(S^2)=0$;
    \item $\tr(S^2\circ r(a))=0$ for all $a \in H$. \qed\\
  \end{enumerate}
\end{theorem}

\begin{proposition}[{\cite[Corollary 2.2]{Ng03}}] \label{p0.1}
  Let $H$ be a finite-dimensional Hopf algebra over $k$ with
  antipode $S$ and $g$ the distinguished group-like
  element of $H$. If\,  $\lcm(o(S^4), o(g))=n$ is an odd integer
  greater than 1, then the subspace
  $$
  H_-=\{u \in H\mid S^{2n}(u)=-u\}
  $$
  has even dimension. \qed\\
\end{proposition}
 The following lemma is useful in our remaining discussion.
\begin{lemma}\label{l0}
  Let $V$ be a finite-dimensional vector space over the field $k$.
  If $T$ is a linear automorphism on $V$ such that $\tr(T)=0$ and
  $o(T)=q^n$ for some prime $q$ and positive integer $n$,
  then
  $$
  q \mid \dim(V)\,.
  $$

\end{lemma}
\begin{proof}
Let $\w \in k$ be a primitive $q^n$th root of unity and
$$
V_i= \{u \in V\mid T(u)=\w^i u\} \quad \mbox{for} \quad i=0, \cdots, q^n-1\,.
$$
Consider the integral polynomial
$$
f(x)=\sum_{i=0}^{q^n-1}\dim(V_i) x^i\,.
$$
Since
$$
0 = \Tr(T) = f(\w)\,,
$$
there exists $g(x) \in \BZ[x]$ such that
$$
f(x)=\Phi_{q^n}(x)g(x)\,
$$
where $\Phi_{q^n}(x)$ is the $q^n$th cyclotomic polynomial.
 Hence,
$$
\dim(V)=f(1)=\Phi_{q^n}(1)g(1)\,.
$$
Since $\Phi_{q^n}(x)=\Phi_q(x^{q^{n-1}})$,
$\Phi_{q^n}(1)=\Phi_q(1)=q$. Thus we have
$$
q \mid \dim(V)\,.
$$
\end{proof}
\section{Unimodularity of Hopf algebras of dimension $2n$}
In this section, we prove that if $H$ is a non-semisimple Hopf algebra over $k$
of dimension $2n$, where $n$ is an odd integer, then $H$ or $H^*$
is not unimodular. This result is  essential to the proof of our
main result in the next section.
\begin{proposition}\label{p1}
  Let $H$ be a finite-dimensional Hopf algebra over the field $k$ with antipode $S$.
  If $H$ is unimodular and $o(S^2)=2$, then
   $$
   4\,\mid \dim(H)\,.
   $$
\end{proposition}
\begin{proof}
Let $\lambda$ be a non-zero right integral of $H^*$. Since $H$ is
unimodular, by \cite[Proposition 2]{Radf94},
$$
\lambda(ab) = \lambda(S^2(b)a)
$$
for all $a, b \in H$. Let
$$
H_i = \{u \in H| S^2(u)=(-1)^i u\} \quad \mbox{for}\quad i=0,1\,.
$$
We claim that $(a,b)=\lambda(ab)$ defines a non-degenerate alternating form on
$H_1$. For any $a, b \in H_1$,
$$
\lambda(ab)=\lambda(S^2(b)a)=-\lambda(ba)\,.
$$
Since $\lambda(u)=\lambda(S^2(u))=-\lambda(u)$ for all $u \in H_1$, $\lambda (H_1)=\{0\}$.
Let $a \in H_1$ such that $\lambda(ab)=0$ for all $b \in
H_1$. Then for all $b \in H_0$, $ab \in H_1$ and so $\lambda(ab)=0$.
By the non-degeneracy of $\lambda$ on $H$, $a=0$. Therefore, $(a,b)=\lambda(ab)$ defines an
non-degenerate alternating bilinear form on $H_1$ and hence $\dim(H_1)$ is even.
Since $o(S^2)=2$, by Theorem \ref{t2}, $\Tr(S^2)=0$ and so $\dim(H_0)=\dim(H_1)$.
Therefore
$$
\dim(H)=\dim(H_0)+\dim(H_1)=2\dim(H_1)
$$
is a multiple of 4.
\end{proof}
\begin{corollary}\label{c1}
  Let $H$ be a Hopf algebra over $k$ of dimension $2n$ where $n$ is an odd integer.
 If $H$ is not semisimple, then $H$ or $H^*$ is not unimodular.
\end{corollary}
\begin{proof}
If both $H$ and $H^*$ are unimodular, by Proposition \ref{p0},
$S^4=id_H$. Since $H$ is not semisimple, by Theorem \ref{t2},
$o(S^2)=2$. It follows from Proposition \ref{p1} that $\dim(H)$ is
then a multiple of 4 which contradicts $\dim(H)=2n$.
\end{proof}
\section{Hopf algebras of dimension $2p$}
In this section, we prove, by contradiction, that non-semisimple Hopf algebras over
$k$ of dimension $2p$, $p$ an odd prime, do not exist. By \cite{Masuoka95}, semisimple
Hopf algebras of dimension $2p$ are
$$
k[\BZ_{2p}],\quad k[D_{2p}]\quad\mbox{and}\quad k[D_{2p}]^*
$$
where $D_{2p}$ is the dihedral group of order $2p$. Our main result completes the
classification of Hopf algebras of dimension $2p$. We begin to
prove our main result with the following lemma.

\begin{lemma}\label{l3}
Let $H$ be a non-semisimple finite-dimensional Hopf algebra over $k$ of dimension $2p$ where $p$ is an odd
prime.  Suppose that $g$ and $\a$ are the distinguished group-like element of $H$ and $H^*$ respectively.
Then
$$
\lcm(o(g), o(\a))=2 \mbox{ or } p\,.
$$
\end{lemma}
\begin{proof}

Since $H$ is not semisimple, by Theorem \ref{t2}, $H^*$ is also not semisimple.
Therefore $|G(H)|$ and $|G(H^*)|$ are strictly less than $2p$. By Nichols-Zoeller theorem \cite{Nich89},
$$
|G(H)|,\,\, |G(H^*)| \in \{1, 2, p\}\,.
$$
It follows from \cite[Lemma 5.1]{Ng02} that
$$
\lcm(|G(H)|, |G(H^*)|)= 1, 2 \mbox{ or } p\,.
$$
Since $\lcm(\,o(g), o(\a)\,)$ divides $\lcm(\,|G(H)|, |G(H^*)|\,)$,
we obtain
$$
\lcm(o(g), o(\a))=1, 2 \mbox{ or } p\,.
$$
By Corollary \ref{c1}, $\lcm(o(g), o(\a)) > 1$ and so the result follows.
\end{proof}
\begin{lemma}\label{l4}
  Let $H$ be a finite-dimensional Hopf algebra over $k$
  and $a \in G(H)$ of order $d$.
  Let $\w \in k$ be a primitive $d$th root of unity and
  $$
  e_i = \frac{1}{d}\sum_{j=0}^{d-1} \w^{-ij}a^j \quad(i=0, \dots, d-1).
  $$
  Then  $\dim(H)/d$ is an integer, $\dim (He_i) = \dim(H)/d$
  and $S^2(He_i) =He_i$ for $i=0, \dots, d-1$. In addition, if $H$
  is not semisimple, then
  $$
  \tr(S^2|_{He_i})=0
  $$
   for $i=0, \dots, d-1$\,.
\end{lemma}
\begin{proof}
Let $B=k[a]$. Then $B$ is a Hopf subalgebra of $H$ and
$\dim(B)=d$. By Nichols-Zoeller theorem, $H$ is a free
$B$-module. In particular, $\dim(H)$ is a multiple of $d$
and
$$
H \cong B^{\dim(H)/d}
$$
as right $B$-modules. Note that $e_0, \dots, e_{d-1}$ are
orthogonal idempotents of $B$ such that
$$
1=e_0+ \dots + e_{d-1}\,,
$$
and $Be_i=ke_i$\,.
Therefore,
$$
He_i \cong B^{\dim(H)/d}e_i =
(Be_i)^{\dim(H)/d}=(ke_i)^{\dim(H)/d}
$$
and so $\dim(He_i)=\dim(H)/d$ for $i=0, \dots, d-1$. Since
$S^2(a)=a$, $S^2(e_i)=e_i$ for $i=0, \dots, d-1$. Therefore,
$$
S^2(He_i) = HS^2(e_i) =He_i\,.
$$
If, in addition, $H$ is not semisimple, by Theorem \ref{t2},
$$
\tr(S^2|_{He_i})=\tr(S^2 \circ r(e_i))=0\,.
$$
for $i=0, \dots, d-1$.
\end{proof}

\begin{theorem}\label{main}
If $p$ is an odd prime, then any Hopf algebra of dimension $2p$ over the field $k$ is semisimple.
\end{theorem}
\begin{proof}
Suppose there exists a non-semisimple Hopf algebra $H$ of
dimension $2p$. Let $g$ and $\a$ be the distinguished group-like
elements of $H$ and $H^*$ respectively. By Corollary \ref{c1}, $g$
and $\a$ can not be both  trivial. By Theorem \ref{t2}, we may
simply assume that $g$ is not trivial and $o(g)=d$. Let $\w \in
k$ be a primitive $d$th of unity and
$$
e_i= \frac{1}{d}\sum_{j=0}^{d-1} \w^{-ij}g^j \quad(i=0, \dots,
d-1).
$$

By Lemma \ref{l3},
$$
\lcm(o(g), o(\a))=2 \mbox{ or } p\,.
$$
\\
If $\lcm(o(g), o(\a))=2$, then $d=2$ and $S^8=id_H$ by Proposition
\ref{p0}. It follows from Theorem \ref{t2} that
$$
o(S^2)=2\quad\mbox{or}\quad 4\,.
$$
It follows from Lemma \ref{l4} that
$$
\dim(He_i)=p \quad \mbox{and}\quad \tr(S^2|_{He_i})=0 \quad(i=0,1)\,.
$$
Since $o(S^2)=2$ or $4$, there exists $j \in \{0, 1\}$ such that
$$
o(S^2|_{He_j})=2 \quad\mbox{or}\quad 4\,.
$$
By Lemma \ref{l0}, $\dim(He_j)$ is even which contradicts that $\dim(He_j)=p$.\\

If $\lcm(o(g), o(\a))=p$, then $d=p$ and $S^{4p}=id_H$. By
Proposition \ref{p0.1}, the subspace
$$
H_-=\{u \in H|S^{2p}(u)=-u\}
$$
has even dimension. On the other hand,
by Lemma \ref{l4}, we have
$$
\dim(He_i)=2 \quad \mbox{and}\quad \tr(S^2|_{He_i})=0 \quad
\mbox{for} \quad i=0,\dots, p-1\,.
$$
Thus, for $i \in \{0, \dots, p-1\}$, there is a basis $\{u_{i}^+, u_{i}^-\}$ for
$He_i$ such that
$$
S^2(u_i^{\pm})=\pm \zeta_i
$$
for some $p$th root of unity $\zeta_i$.  Thus, $\{u_0^-, u_1^-,
\dots, u_{p-1}^-\}$ forms a basis of $H_-$ and so
$$
\dim(H_-)=p\,,
$$
a contradiction!
\end{proof}
\vspace{0.5cm}
\begin{center}
 {\bf Acknowledgement}
\end{center}
The author would like to thank L. Long for her useful suggestions
to this paper.
\bibliographystyle{amsalpha}
\providecommand{\bysame}{\leavevmode\hbox to3em{\hrulefill}\thinspace}
\providecommand{\MR}{\relax\ifhmode\unskip\space\fi MR }
\providecommand{\MRhref}[2]{%
  \href{http://www.ams.org/mathscinet-getitem?mr=#1}{#2}
}
\providecommand{\href}[2]{#2}

\end{document}